\documentclass[journal]{IEEEtran}

\ifCLASSINFOpdf
\else
   \usepackage[dvips]{graphicx}
\fi
\usepackage{url}

\usepackage{amsmath,amsfonts,amssymb,amsthm}

\usepackage{subcaption}

\usepackage{graphicx}

\begin{document}

\title{On Efficient Sampling Schemes for the Eigenvalues of Complex Wishart Matrices}

\author{Peter J. Forrester
\thanks{Submitted 18th December 2023. This research is part of the program of study supported
	by the Australian Research Council Discovery Project grant DP210102887. }
\thanks{Peter J. Forrester  is with the School of Mathematics and Statistics, University of Melbourne, Victoria 3010, Australia.(e-mail: pjforr@unimelb.edu.au).}}

\markboth{IEEE Signal Processing Letters}
{Shell \MakeLowercase{\textit{et al.}}: Bare Demo of IEEEtran.cls for IEEE Journals}
\maketitle

\begin{abstract}
The paper ``An efficient sampling scheme for the eigenvalues of dual Wishart matrices'', by I.~Santamar\'ia and V.~Elvira, 
[\emph{IEEE Signal Processing Letters}, vol.~28, pp.~2177--2181, 2021] \cite{SE21}, poses the question of efficient sampling from the 
eigenvalue probability density function of the $n \times n$ central
complex Wishart matrices with variance matrix equal to the identity. Underlying such complex Wishart matrices is a rectangular
$R \times n$ $(R \ge n)$ standard complex Gaussian matrix, requiring then $2Rn$ real random variables for their
generation. The main result of \cite{SE21} gives a formula involving just two classical
distributions specifying the two eigenvalues in the case $n=2$. The purpose of this Letter is to point out that existing results in
the literature give two distinct ways to efficiently sample the eigenvalues in the general $n$ case. One is in terms of the eigenvalues of
a tridiagonal matrix which factors as the product of a bidiagonal matrix and its transpose, with the  $2n+1$ nonzero entries of the latter
given by (the square root of) certain chi-squared random variables. The other is as the generalised eigenvalues for a pair of bidiagonal
matrices, also containing a total of $2n+1$ chi-squared random variables. Moreover, these characterisation persist in the case
of that the variance matrix consists of a single spike, and for the case of real Wishart matrices.
\end{abstract}


\IEEEpeerreviewmaketitle

\section{Introduction}

\IEEEPARstart{T}{his} Letter is a follow up to the work \cite{SE21} by Santamar\'ia and Elvira relating  to the efficient sampling of the the eigenvalues of
central complex Wishart matrices with variance matrix equal to the identity. 
Random matrices of this type are of the form $\mathbf W = \mathbf H^{\rm H} \mathbf H$, where $ \mathbf H$ is an $R \times n$ ($R \ge n$) standard complex 
Gaussian matrix (i.e.~each entry is independent and distributed as a standard complex normal) and $ \mathbf H^{\rm H}$ denotes
the Hermitian conjugate of $ \mathbf H$. The matrix $\mathbf H$ may be considered as a data matrix of $R$ measurements of each of $n$ quantities forming the columns,
where each quantity is distributed as a standard complex normal.
Extending the classical work of Wishart \cite{Wi28} in the real case, Goodman \cite{Go63}
derived that the element probability density function (PDF) of $ \mathbf W$ is proportional to $(\det  \mathbf W)^{R-n} e^{- {\rm Tr} \,  \mathbf W}$, and moreover
that the PDF of the ordered eigenvalues $\lambda_1 > \lambda_2 > \cdots > \lambda_n > 0$ is equal to
\begin{equation}\label{(1)}
{ n! e^{ - \sum_{j=1}^n \lambda_j} \prod_{l=1}^n \lambda_l^{R - n} \over \prod_{j=2}^n \Gamma(1 + j) \Gamma(j+R - n)}
\prod_{1 \le j < k \le n} (\lambda_k - \lambda_j)^2.
\end{equation}
Motivated by the relevance of (\ref{(1)}) to MIMO wireless communications, the authors of  \cite{SE21} posed the question of efficient
methods for its sampling.

Rather than constructing $\mathbf W$ directly from its definition in terms of the standard complex Gaussian matrix $\mathbf H$ which requires
$Rn$ complex random variables for its generation, 
in \cite{SE21} the Gram-Schimdt decomposition was applied to the columns of $\mathbf H$ to give the matrix decomposition
$ \mathbf H =  \mathbf Q \mathbf R$. Here $ \mathbf Q$ is an $R \times n$ matrix with columns forming an orthonormal basis of the column space of $ \mathbf H$, while $ \mathbf R$
is an upper triangular $n \times n$ matrix with positive real diagonal entries. Multiplying the Hermitian conjugate of $\mathbf H$ times
$\mathbf H$, one sees that the dependence on $\mathbf Q$ cancels,  showing that $ \mathbf W =  \mathbf R^{\rm H}  \mathbf R$. Thus now
$ \mathbf W$ is expressed as the product of an $n \times n$ matrix and its Hermitian conjugate, rather than an $R \times n$ matrix (on the right)
times its Hermitian conjugate. Most importantly, a result due to Bartlett in the real case \cite{Ba33}, and extended to the complex
case by Goodman \cite[part of the proof of Th.~5.1]{Go63} gives that the square of the $j$-th diagonal entry of $ \mathbf R$ has gamma distribution
${\rm Gamma}(R-j+1,1)$, or equivalently chi-squared distribution ${1 \over 2} \chi_{2(R-j+1)}^2$. Also, all strictly upper triangular triangular entries
are standard complex normal random variables, and are thus unchanged in distribution from the entries of $ \mathbf H$.

Write $ \mathbf R = [r_{ij} ]_{i,j=1}^n$, with $r_{i,j} =0$ for $i > j$. In the case $N = 2$ the above theory gives \cite[Eq.~(3)]{SE21}
\begin{equation}\label{(3)}
W = \begin{bmatrix} r_{11}^2 & r_{11} r_{12} \\
r_{12}^* r_{11} & |r_{12}|^2 + r_{22}^2 \end{bmatrix},
\end{equation}
where the $r_{ij}$ have distribution $r_{11}^2 \mathop{=}\limits^{\rm d} {1 \over 2} \chi_{2R}^2$, $r_{22}^2 \mathop{=}\limits^{\rm d} {1 \over 2} \chi_{2(R-1)}^2$ and
$r_{12}  \mathop{=}\limits^{\rm d}  {\mathcal CN}(0,1)$. The characteristic polynomial for (\ref{(3)}) is \cite[Eq.~(4)]{SE21}
\begin{align}\label{(4)}
\det (\mathbf W - \lambda \mathbf I_2) & = \lambda^2 - \lambda {\rm tr} (\mathbf W) + \det (\mathbf W) \nonumber \\
& = \lambda^2 - \lambda (| r_{12} |^2 + r_{11}^2 + r_{22}^2) + r_{11}^2 r_{22}^2 .
\end{align}
As noted in \cite[below Eq.~(6)]{SE21}, it follows from (\ref{(4)}) together with the facts that $r_{11}^2, r_{22}^2$  are chi-squared distributed, as is
$| r_{12} |^2  \mathop{=}\limits^{\rm d}  {1 \over 2} \chi_{2}^2$
that $ {\rm tr} (\mathbf W)  \mathop{=}\limits^{\rm d} {1 \over 2} \chi_{4R}^2$. Less immediate but nonetheless true
is the fact that with $\eta := \det(\mathbf W)/({1 \over 2} {\rm tr} (\mathbf W))^2$ and $s:= 1 - \eta$, one has that \cite[Eq.~(6)]{SE21}
\begin{equation}\label{(5)}
s  \mathop{=}\limits^{\rm d}  {\rm Beta}(3/2,R-1);
\end{equation}
see the following paragraph.

The eigenvalues $\lambda_1 > \lambda_2 > 0$ are the two zeros of the quadratic (\ref{(4)}), suitably ordered and given in terms of $s$ as defined above and $t:=\lambda_1 + \lambda_2 =  {\rm tr} (\mathbf W) $ by
\begin{equation}\label{(6)}
\lambda_1 = {1 \over 2} t (1 + \sqrt{s}), \qquad \lambda_2 = {1 \over 2} t (1 - \sqrt{s}).
\end{equation}
Moreover, changing variables in (\ref{(1)}) to $t$ and $s$ shows that $t$ and $s$ are independently distributed, with distributions as already noted above (as must be). Consequently (\ref{(6)})
provides an efficient way to sample from the $n=2$ eigenvalue PDF (\ref{(1)}), which is the main result from \cite{SE21}. Furthermore \cite{SE21} notes that (\ref{(6)}) applies equally as
well to the $n=2$ eigenvalue PDF for real Wishart matrices (see e.g.~\cite[Eq.~(3.16) with $\beta = 1$, $m=2$ and $n=R$]{Fo10}, now with $t   \mathop{=}\limits^{\rm d}  \chi_{2R}^2$ and
$s \mathop{=}\limits^{\rm d}  {\rm Beta}(1,(R-1)/2)$.

\section{The case of general $n$}

In motivating the work \cite{SE21},  one reads in the Introduction therein (below  Eq.~(1)) that ``to the best of the best of the authors' knowledge there are no efficient procedures to sample from [the $n=2$ case of our
(\ref{(1)})]''. In the conclusion section of \cite{SE21} it is stated that ``future work will consider  extension to other cases, e.g.~non-central Wishart distributions or Wishart matrices of dimension greater than two''. It is the purpose of this Letter to point out  that there is a known efficient random matrix  formulation to sample from (\ref{(1)}) in the general $n$ case, and also for the single spiked Wishart extension of (\ref{(1)}).

We consider first the problem of efficiently sampling from (\ref{(1)}) for general $n \ge 2$. This is available in the literature from the 2002 work of Dumitriu and Edelman \cite{DE02}. Starting with a complex ($\beta = 2$) or real ($\beta = 1$) Wishart matrix $\mathbf W$, the main idea from \cite{DE02} (see also \cite{Si85} in the real case) is that the underlying Gaussian matrices $\mathbf R$ can be reduced, using Householder transformations, to an $n \times n$ bidiagonal matrix $\mathbf B$, such that $\mathbf B^{\rm H} \mathbf B$ has the same joint eigenvalue distribution of $\mathbf W$. Remarkably, the nonzero elements in $\mathbf B$ are all independently distributed as (the square root of) certain chi-squared random variables, which permit for a generalisation in their degrees of freedom allowing for $\beta$ (originally labelling the number field) to be a continuous positive real parameter. Explicitly, the  matrix $\mathbf B = \mathbf B_\beta$ is given by
 \begin{equation}\label{Bd}
 \mathbf B_\beta =
{1 \over \sqrt{\beta} } \begin{bmatrix} \chi_{\beta R} & & &  \\
\chi_{\beta (n - 1)} & \chi_{\beta (R - 1)} & &  \\
 & \chi_{\beta (n - 2)} & \chi_{\beta (R - 2)} &  & \\
 & \ddots & \ddots &  & \\
& & \chi_\beta & \chi_{\beta (R - n + 1)}
\end{bmatrix}.
\end{equation}
The main theorem from \cite{DE02} is that $\mathbf B_\beta^{\rm H}  \mathbf B_\beta$ has eigenvalue PDF generalising (\ref{(1)}), specified by the functional form
\begin{equation}\label{1X}
{ 1 \over \mathcal N_{n,R,\beta}} e^{-\beta \sum_{j=1}^n  \lambda_j/2 } \prod_{l=1}^n \lambda_l^{\beta a /2}
\prod_{1 \le j < k \le n} |\lambda_k - \lambda_j|^\beta.
\end{equation}
Here $a := R - n + 1 - 2/\beta$ (with $R \ge n$ also permitted to be continuous) and $ \mathcal N_{n,R,\beta}$ denotes the normalisation (which is known in terms of a product
of gamma functions, being a limiting case of the Selberg integral; see e.g.~\cite[\S 4.7]{Fo10}). Hence sampling from (\ref{1X}) can be carried out by computing the
eigenvalues of the $n \times n$ tridiagonal matrix obtained by computing the transpose times matrix product of the bidiagonal matrix (\ref{Bd}), with the latter
consisting of (the square root of) $2n-1$ independent chi-squared random variables. With $\beta = 2$, this answers the question of efficiently sampling from (\ref{(1)}) for general 
$n \ge 2$ as posed
in the  conclusion section of \cite{SE21}. Recall that the original Gaussian matrix construction of the complex Wishart matrix involves the at least quadratic in $n$, $R  n$ ($R \ge n$), independent complex 
Gaussian entries to obtain the matrix $\mathbf H$, then forming $\mathbf W = \mathbf H^{\rm H} \mathbf H$ and computing the eigenvalues of the dense matrix $\mathbf W$.

The fact that $\mathbf B^{\rm H}_\beta \mathbf B_\beta := \mathbf T_\beta$ is tridiagonal implies that the corresponding characteristic polynomial satisfies a three term recurrence \cite{DE02}.
In this regard, introduce the notation
$$
\mathbf B_\beta :=
\left [ \begin{array}{cccc}
x_R & & & \\
y_{n-1} & x_{R-1} & & \\
\ddots & \ddots & & \\
& y_1 & x_{R-n+1} & \end{array} \right ],
$$
$$
\mathbf T_\beta :=
\left [ \begin{array}{cccc}
a_n &b_{n-1} & & \\
b_{n-1} & a_{n-1} &b_{n-2} & \\
\ddots & \ddots &\ddots & \\
& b_2 & a_{2} & b_1 \\
& & b_1 & a_1 \end{array} \right ],
$$
so that
\begin{equation}\label{2.67'}
a_n = x_R^2, \quad a_i = y_i^2 + x_{R-n+i}^2, \quad
b_i = y_i x_{R-n+i+1},
\end{equation}
with $i=n-1,n-2,\dots,1$.
Denote by $P_k(\lambda)$ the characteristic polynomial of the bottom right $k \times k$ submatrix of $ \mathbf T_\beta$. By expanding the
determinant defining $P_k(\lambda)$ along its first row, one obtains the three-term recurrence
\begin{equation}\label{2.68}
P_k(\lambda) = (\lambda - a_k) P_{k-1}(\lambda) - b_{k-1}^2  P_{k-2}(\lambda),
\end{equation}
which with the initial conditions $P_{-1}(\lambda) = 0$ and $P_0(\lambda) = 1$  allows for the recursive computation of the characteristic polynomial $P_n(\lambda)$ of
$B_\beta$.

Suppose now that in the construction of $\mathbf W $ we replace the complex (or real) $R \times n$ Gaussian matrix $\mathbf H$ by $\mathbf H \mathbf \Sigma^{1/2}$, where $\mathbf  \Sigma$ is
a fixed positive definite matrix. Denote the eigenvalues of $\mathbf  \Sigma$ by $\{ \sigma_j^2 \}_{j=1,\dots,n }$. From the fact that $\mathbf H$ consists of independent Gaussian entries, it follows that
$(\mathbf H \mathbf \Sigma^{1/2})^{\rm H} \mathbf H \mathbf \Sigma^{1/2}$ is unitary equivalent and thus has the same eigenvalues as the same matrix product with $\mathbf  \Sigma$ replaced by the
diagonal matrix of its eigenvalues. In this latter setting, we have that the variance of the Gaussian elements in column $j$ is, upon multiplication by ${\rm diag}( \mathbf \Sigma^{1/2})$, equal to $\sigma_j^2$,
so being considered is the covariance matrix of a data matrix in which column $j$ has entries which are  zero mean complex normals with variance $\sigma_j^2$. 
The case where only one of the variances (say $\sigma_1$) is different from unity is referred to as the multiplicative rank 1 case (or as a single spike); generally rank 1 perturbations hold a special place in random
matrix theory as reviewed recently in \cite{Fo22}. In the present setting we know from \cite{Fo13} that the Householder reduction strategy still applies, with the only effect being that the top left
entry in (\ref{Bd}) is now to be multiplied by $\sigma_1$. With this adjustment, one has in particular that the corresponding characteristic polynomial can be computed via the (random) three term
recurrence (\ref{2.68}).

It turns out that there is a distinct random three term recurrence giving the characteristic polynomial  in the above specified multiplicative rank 1 case \cite{FR02b,Fo13}. This comes about by noting
that in the complex or real cases the Gaussian matrix product defining $\mathbf W$ permits the decomposition
\begin{equation}\label{2.69}
(\mathbf H \mathbf \Sigma^{1/2})^{\rm H} \mathbf H \mathbf \Sigma^{1/2} = \tilde{\mathbf H}^{\rm H} \tilde{\mathbf H} + \sigma_1^2 \mathbf x \mathbf x^{\rm H},
\end{equation}
where $\tilde{\mathbf H}$ is the $(R -1 ) \times n$ matrix obtained from $\mathbf H$ by deleting its first column.  The column vector $\mathbf x$ has $n$ components which are independent
complex or real standard normals as appropriate.  A key observation from \cite{FR02b} (see also \cite{DEKV13}) is that the RHS of (\ref{2.69}) is unitary equivalent to
\begin{equation}\label{2.70}
{\rm diag}( \tilde{\mathbf H}^{\rm H} \tilde{\mathbf H}) +  \sigma_1^2 \mathbf x \mathbf x^{\rm H}
\end{equation}
and moreover the eigenvalues of this matrix depend on the entries $x_j$ of $\mathbf x$ only through the squared modulus $|x_j|^2  \mathop{=}\limits^{\rm d} {1 \over \beta} \chi_{\beta}^2$
($\beta = 1$ real case and $\beta = 2$ complex case). Thus the eigenvalue problem is well defined from general $\beta > 0$, in which case we know that the PDF of the eigenvalues of
$\tilde{\mathbf H}^{\rm H} \tilde{\mathbf H}$ is given by (\ref{Bd}) with the parameter $a$ reduced by 1 to account for the size of $\tilde{\mathbf H}$ being $(R -1 ) \times n$. In \cite{FR02b}
an explicit formula for the joint PDF of the eigenvalues of $\tilde{\mathbf H}^{\rm H} \tilde{\mathbf H}$ and of $(\mathbf H \mathbf \Sigma^{1/2})^{\rm H} \mathbf H \mathbf \Sigma^{1/2}$ along
with a three term recurrence which gives the characteristic polynomials of each. In relation to the latter, introduce the random variables
$$
\begin{array}{ll}
&\tilde{a}_j   \mathop{=}\limits^{\rm d} {1 \over \beta} \chi_{(R+1-j)\beta/2}^2 \: \: (j=1,\dots,n-1),  \\
&\tilde{a}_n   \mathop{=}\limits^{\rm d} {\sigma_1^2 \over \beta} \chi_{ (R+1-n)\beta/2}^2  \\
&\tilde{b}_j   \mathop{=}\limits^{\rm d} {1 \over \beta} \chi_{j\beta/2}^2 \: \: (j=1,\dots,n-2),  \\
& \tilde{b}_{n-1}   \mathop{=}\limits^{\rm d} {\sigma_1^2 \over \beta} \chi_{ (n-1)\beta/2}^2 .
\end{array}
$$
Then the results of \cite[\S 5.2]{FR02b} tell us that the random monic polynomials defined by 
 \begin{equation}\label{2.71}
 B_j(x) = (x - \tilde{a}_j )  B_{j-1}(x)  - \tilde{b}_{j-1}  B_{j-2}(x) \quad (j=1,\dots,n),
 \end{equation}
 with $B_{-1}(x) = 0$ and $B_0(x) = 1$
 is such that $B_{n-1}(x)$ and $B_{n}(x)$ are the sought characteristic polynomials. In \cite{Fo13} it is noted that in general (\ref{2.71}) is satisfied by 
 $ B_j(x) = \det ( x M_j - L_j)$ where $L_j$ and $M_j$ are the top $j \times j$ blocks of the bidiagonal matrices
 $$
\mathbf L:= \begin{bmatrix}\tilde{a}_1 & 1 & & &\\
 & \tilde{a}_2 & 1 & &\\
& & \ddots & \ddots  &\\
& & & \tilde{a}_{n-1} & 1 \\
& & &  & \tilde{a}_n \end{bmatrix}, 
$$
$$
\mathbf  M:= \begin{bmatrix}1 &  & & &\\
-\tilde{b}_1 & 1 &  & &\\
 & -\tilde{b}_2 & 1 && \\
 &  & \ddots & \ddots  &\\
& & &  -\tilde{b}_{n-1} & 1 \end{bmatrix}.
 $$
 Hence from this approach generation of the PDF of $\{\lambda_j\}$ as specified by (\ref{Bd})  is reduced to a generalised eigenvalue problem for
 a pair of bidiagonal matrices.

 \section{Simulation of marginals and comparison with theory}
 
 \begin{figure}[h]
\begin{subfigure}[b]{0.90\columnwidth}
\centerline{\includegraphics[width=\columnwidth]{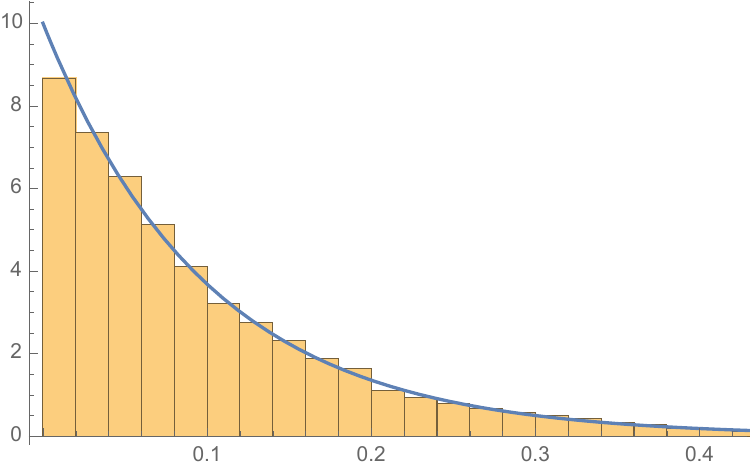}}
\end{subfigure}

\begin{subfigure}[b]{0.90\columnwidth}
\centerline{\includegraphics[width=\columnwidth]{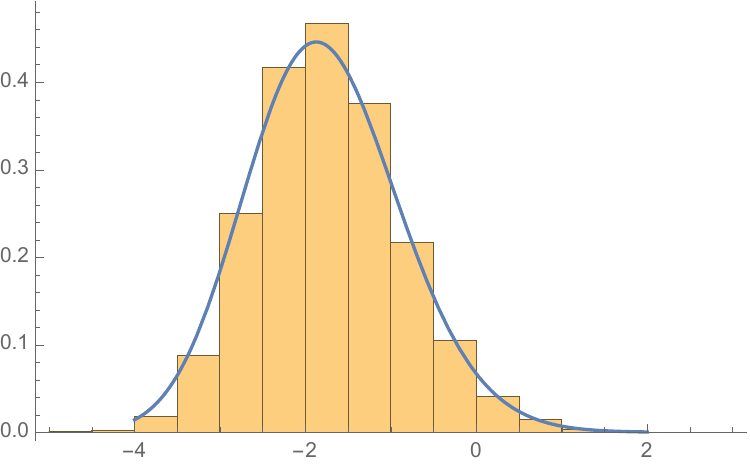}}

\end{subfigure}
\caption{Empirical PDFs generated by sampling $2 \times 10^4$ times the smallest and largest eigenvalue respectively from the random tridiagonal matrix $\mathbf T_\beta $ with parameters $n=R=30$ and $\beta = 2$,    compared against theoretical predictions. }
\label{EF.1}
\end{figure}

 The functional form of the marginal distribution of a single ordered eigenvalue as implied by (\ref{(1)}), in particular the smallest and largest, has been the subject of much
 attention over the years. The interest comes from both applications to wireless communications, and to quantum transport; see e.g.~the Introduction of 
 \cite{ZCW09} or \cite{FK19}. The reference \cite{FK19} takes as its objective of expressing the the PDF of the largest eigenvalue as implied by (\ref{(1)}) with  in a terminating form suitable for
 further analysis, made
 possible by identifying appropriate basis functions. The simplest result of this sort was identified long ago in relation to the PDF of smallest eigenvalue,
 $p_{\rm min}(x;n,R)$ say \cite{FH94}. In particular, one has from \cite{FH94} that
  \begin{equation}\label{pr}
 p_{\rm min}(x;n,n) = n e^{-x n}.
 \end{equation}
 One notes from from (\ref{pr}) the scaling limit $\lim_{n \to \infty} {1 \over n}  p_{\rm min}(x;n,n) = e^{-x}$. However, while such a limit formula is well known
 to exist for $p_{\rm min}(x;n,R)$ with $a:=R-n$ fixed, its functional form depends on $a$   \cite{FH94}. On the other hand, considering instead the largest eigenvalue, and
 recentering and rescaling according to 
  \begin{equation}\label{rr}
  \lambda_1 \mapsto (\lambda_1 - 4 N - 2a)/2 (2 N)^{1/3},
  \end{equation} 
one has that with $a$ fixed the PDF is independent of $a$, and limits to the so-called Tracy-Widom $\beta = 2$ PDF at a rate which is proportional to $1/N^{2/3}$;
on this last point see \cite{FT18}.

The tridiagonal matrix $\mathbf T_\beta |_{\beta = 2} = \mathbf B^{\rm H}_\beta \mathbf B_\beta |_{\beta = 2}$ with $ \mathbf B_\beta$ the bidiagonal matrix (\ref{Bd}), is well suited
to efficiently generating multiple samples of specific ordered eigenvalues from (\ref{(1)}). Using the software package Mathematica, this is further facilitated by the {\tt Eigenvalue} command
for computing only the largest or the smallest eigenvalue. Also, the functional form of the Tracy-Widom $\beta = 2$ distribution is inbuilt \cite{WMath}. Our Figure \ref{EF.1} reports in
graphical form normalised histograms for the PDF of the smallest and (recentered and rescaled) largest eigenvalue generated from $\mathbf T_\beta |_{\beta = 2}$ with $R=n=30$, using $2 \times 10^4$
samples. Superimposed are the theoretical result (\ref{pr}) in the case of the smallest eigenvalue, and the Tracy-Widom $\beta = 2$ PDF, as applies to the limiting recentered and rescaled
largest eigenvalue.

\section*{Acknowledgment}
Correspondence relating to an earlier draft of this work by I.~Santamar\'ia is appreciated.

\end{document}